\documentclass[12pt, a4paper]{article}
\usepackage[latin1]{inputenc}
\usepackage[english]{babel}
\usepackage{indentfirst}
\usepackage{amstext}
\usepackage{setspace}
\usepackage{amsfonts}
\usepackage{textcomp}
\usepackage{amssymb}
\usepackage{amscd}
\usepackage{epsf}
\usepackage{graphicx}
\usepackage{epsfig}

\newcommand {\demo}{\hskip -0.6cm{\bf Proof:  }}
\newcommand {\fim}{\hfill{$\square$}\vskip 1pc}

\newcommand {\N}{\mathbb{N}}

\newcommand {\Lh}{\mathcal{L}}

\newtheorem{theorem}{Theorem}[section]
\newtheorem{lema}[theorem]{Lemma}
\newtheorem{corolario}[theorem]{Corollary}
\newtheorem{definition}[theorem]{Definition}
\newtheorem{proposition}[theorem]{Proposition}
\newtheorem{example}[theorem]{Example}

\begin{document}

\title{Perron-Frobenius operators and representations of the Cuntz-Krieger algebras for infinite matrices}
\maketitle
\begin{center}
{\large Daniel Gonçalves\footnote{The author would like to thank Wael Bahsoun for bringing Kawamura´s, \cite{kawamura}, paper to his attention.} and Danilo Royer \footnote{partially supported by a grant from FAPESC, Brazil}}\\
\end{center}  
\vspace{8mm}

\abstract In this paper we extend work of Kawamura, see \cite{kawamura}, for Cuntz-Krieger algebras $O_A$ for infinite matrices $A$. We generalize the definition of branching systems, prove their existence for any given matrix $A$ and show how they induce some very concrete representations of $O_A$. We use these representations to describe the Perron-Frobenius operator, associated to an nonsingular transformation, as an infinite sum and under some hypothesis we find a matrix representation for the operator. We finish the paper with a few examples.
\onehalfspace

\section{Introduction}

The interactions between the theory of dynamical systems and operator algebras are one of the main venues in modern mathematics. Exploring this interplay Kawamura, see \cite{kawamura}, recently showed that the theory of representations of the Cuntz-Krieger algebras is closely related to the theory involving the Perron-Frobenius operator. The work of Kawamura is done for the Cuntz-Krieger algebras $O_A$, for finite matrices $A$. In this paper we generalize many of the results in \cite{kawamura} for the Cuntz-Krieger algebras for infinite matrices (a concept introduced by Exel and Laca in \cite{exellaca}). For example, under some mild assumptions, we are able to give a explicit characterization of the Perron-Frobenius operator, associated to a nonsingular transformation, as a infinite sum, using a representation of an infinite Cuntz-Krieger algebra. In our efforts to generalize the notions of \cite{kawamura} we found two problems with the work done in there that we believe are worth mentioning. First is the necessity of an extra hypothesis in the definition of a branching function system given in \cite{kawamura}. The other problem is in the statement of theorem 1.2 of \cite{kawamura}, where $BA$ should read $A^TB$. We will deal with both these cases when introducing our generalized versions of the theory of \cite{kawamura}.

We organize the paper in the following way: In the remaining of the introduction we quickly recall the reader the main definitions of \cite{kawamura} and show the need for an extra hypothesis in the definition of a branching function system. In section 2, we define branching systems for infinite matrices $A$, which we denote by $A_\infty$. We deal with the existence of $A_\infty$ branching systems for any given matrix $A$ (infinite or not) and show how they induce representations of $O_A$ in section 3. Next, in section 4, we use the representations introduced in section 3 to describe the Perron-Frobenious operator as an infinite sum; we also present the generalizide and corrected version of theorem 1.2 of \cite{kawamura} in this section. We finish the paper in section 5 with a few examples.

Given a measure space $(X,\mu)$, let $L_p(X,\mu)$ be the set of all complex valued measurable functions $f$ such that $\|f\|_p<\infty$. For a nonsingular transformation $F:X\rightarrow X$ (that is, $\mu(F^{-1}(A))=0$ if $\mu(A)=0$) let $P_F:L_1(X,\mu)\rightarrow L_1(X,\mu)$ be the Perron-Frobenius operator, that is, $P_F$ is such that $$\int\limits_AP_F\psi(x)d\mu=\int\limits_{F^{-1}(A)}\psi(x)d\mu$$ for each measurable subset $A$ of $X$, for all $\psi \in L_1(X,\mu)$. Notice that, for $\psi \in L_1(X,\mu)$, $P_F(\psi)$ is the Radon-Nikodym derivative of the measure $\mu_{P_F}$, given by $\mu_{P_F}(A)=\int\limits_{F^{-1}(A)}\psi(x)d\mu$, with respect to $\mu$ (see \cite{lasota} for more details about the Perron-Frobenius operator).

In order to describe the Perron-Frobenius operators and representations of the Cuntz-Krieger algebras, Kawamura, in \cite{kawamura}, introduces A-branching function system on a measure space $(X,\mu)$: a family $(\{f_i\}_{i=1}^N, \{D_i\}_{i=1}^N)$ of measurable maps and measurable subsets of $X$ respectively, together with a nonsingular transformation $F:X\rightarrow X$ such that $f_i:D_i\rightarrow f_i(D_i)=R_i$, $\mu(X\setminus \bigcup\limits_{i=1}^N R_i)=0$, $\mu(R_i\cap R_j)=0$ for all $i\neq j$, there exists the Radon-Nikodym derivative $\Phi_{f_i}$ of $\mu\circ f_i$ with respect to $\mu$ and $\Phi_{f_i}>0$ almost everywhere in $D_i$ for $i-1,..,N$, $F\circ f_i=id_{D_i}$ for each $i\in \N$ and $\mu(D_i\setminus\bigcup\limits_{j:a_{ij}=1} R_i) = 0$, where $a_{ij}$ are the entries of the matrix A defining $O_A$. 

Next, a family $\{S(f_i)\}_{i=1}^N$ of partial isometries in $L_2(X,\mu)$ is defined by $S(f_i)(\phi)=\chi_{R_i}\cdot(\Phi_F)^{\frac{1}{2}}\cdot\phi\circ F$, where $\chi_{R_i}$ denotes the characteristic function of $R_i$, and a representation of $O_A$ in $L_2(X,\mu)$ is obtained by defining $\pi_f(s_i)=S(f_i) \ \ (i=1..N)$, (where $s_i$ is one of the generating partial isometry in $O_A$), and using the universal property of $O_A$. But it happens that the definition given above for an A-branching function system is not enough to guarantee that we get a representation of $O_A$, in fact, it is not enough to prove most of the theorems in \cite{kawamura}. For example, let $X=[0,2]$, $\mu$ be the Lebesgue measure, $R_1=[0,1]=D_1, R_2=[1,2]=D_2$, $F:X\rightarrow X$ defined by $F(x)=x$ for each $x\in [0,2]$ (so, $f_i(x)=x$ for each $x\in D_i$) and $A=\left(\begin{array}{cc}1 & 1\\1 & 0\end{array}\right)$. Following $\cite{kawamura}$, $(\{f_i\}_{i=1}^2, \{D_i\}_{i=1}^2)$ is an A-branching function system, but $S(f_1)^*S(f_1)(\phi)=\chi_{[0,1]}\cdot \phi$ and $(S(f_1)S(f_1)^*+S(f_2)S(f_2)^*)(\phi)=\chi_{[0,2]}\cdot \phi$, for each $\phi\in L_2(X,\mu)$, so that $S(f_1)^*S(f_1)\neq\sum\limits_{i=1}^2S(f_i)S(f_i)^*$. Therefore, the existence of a representation of $O_A$ in $L_2([0,2],\mu)$ is not guaranteed.

As we seen, we need to add some extra hypothesis to the definition of an A-branching function system. Namely, we also have to ask that $\mu(\bigcup\limits_{j:a_{ij}=1}R_j\setminus D_i)=0$, for each $i=1,..,N$. We should mention that this extra condition is satisfied in all the examples given in \cite{kawamura}. With this new definition of an A-branching function system in mind, we are now able to generalize it to the countable infinite case.

\section{$A_\infty$-branching systems}

For a measure space $(X,\mu)$ and for measurable subsets $Y,Z$ of $X$, we write $Y\stackrel{\mu -a.e.}{=}Z$ if $\mu(Y\setminus Z)=0=\mu(Z\setminus Y)$ or equivalently, if there exists $Y',Z' \subset X$ such that $Y\cup Y'=Z\cup Z'$ with $\mu(Y')=0=\mu(Z')$. 

Let $A$ be an infinite matrix, with entries $A(i,j)\in \{0,1\}$, for $(i,j)\in \N\times \N$, and let $(X,\mu)$ be a measurable space. For each pair of finite subsets $U,V$ of $\N$ and $j\in \N$ define $$A(U,V,j)=\prod\limits_{u\in U}A_{uj}\prod\limits_{v\in V}(1-A_{vj}).$$

\begin{definition}\label{brancsystem}
An $A_\infty$-branching system on a $\sigma$-finite measure space $(X,\mu)$ is a family $(\{f_i\}_{i=1}^\infty, \{D_i\}_{i=1}^\infty)$ together with a nonsingular transformation $F:X\rightarrow X$ such that:
\begin{enumerate}

\item $f_i:D_i\rightarrow R_i$ is a measurable map, $D_i,R_i$ are measurable subsets of $X$ and  $f_i(D_i)\stackrel{\mu-a.e.}{=} R_i$ for each $i\in \N$; 
\item $F$ satisfies $F\circ f_i=id_{D_i}$  $\mu-a.e.$ in $D_i$ for each $i\in \N$; 
\item $\mu(R_i\cap R_j)=0$ for all $i\neq j$;
\item $\mu(R_j\cap D_i)=0$ if $A(i,j)=0$ and $\mu(R_j\setminus D_i)=0$ if $A(i,j)=1$;
\item for each pair $U,V$ of finite subsets of $\N$ such that $A(U,V,j)=1$ only for a finite number of $j's$, 

$$\bigcap\limits_{u\in U}D_u\bigcap\bigcap\limits_{v\in V}(X\setminus D_v)\stackrel{\mu-a.e}{=}\bigcup\limits_{j\in \N:A(U,V,j)=1}R_j.$$

\item There exists the Radon-Nikodym derivatives $\Phi_{f_i}$ of $\mu\circ f_i$ with respect to $\mu$ in $D_i$ and $\Phi_{f_i^{-1}}$ of $\mu\circ f_i^{-1}$ with respect to $\mu$ in $R_i$.
\end{enumerate}
\end{definition}

The existence of the Radon-Nikodym derivative $\Phi_{f_i}$ of $\mu\circ f_i$ with respect to $\mu$ in $D_i$ together with the fact that $F\circ f_i=Id_{D_i}$ $\mu-a.e.$ imply that $f_i\circ F_{|_{R_i}}=Id_{R_i}$ $\mu-a.e$.
So, the function $f_i$ is $\mu-a.e.$ invertible, with inverse $f_i^{-1}:=F_{|_{R_i}}$. These are the functions that appear in condition $6$ above. If follows from the same condition that $\Phi_{f_i}$ and $\Phi_{f_i^{-1}}$ are measurable functions in $D_i$ and $R_i$ respectively. We will also consider these functions as measurable functions in $X$, defining it as being zero out of $D_i$ and $R_i$, respectively. 

The functions $\Phi_{f_i}$ and $\Phi_{f_i^{-1}}$ are nonnegative $\mu$-a.e., because $\mu$ is a (positive) measure. It is possible to show (by using the following proposition) that $\Phi_{f_i}>0$ and $\Phi_{f_i^{-1}}>0$ $\mu-$a.e. in $D_i$ and $R_i$, respectively. The same proposition shows that $\Phi_{f_i}(x)\Phi_{f_i^{-1}}(f_i(x))=1$ $\mu$-almost everywhere in $D_i$. This equality will be used in the next section.

\begin{proposition}\label{prop1}
Let $(X,\mu)$ be a $\sigma$-finite measure space and $Y,Z$ measurable subsets of $(X,\mu)$. Consider two measurable maps $f:Y\rightarrow Z$ and $g:Z\rightarrow X$, and suppose that there exists the Radon-Nikodym derivatives $\Phi_f$ of $\mu\circ f$ with respect to $\mu$ in $Y$ and $\Phi_g$ of $\mu\circ g$ with respect to $\mu$ in $Z$. Suppose also that $\mu\circ f$ and $\mu\circ g$ are $\sigma$-finite. Then there exists the Radon-Nikodym derivative $\Phi_{g\circ f}$ of $\mu\circ(g\circ f)$ with respect to $\mu$ in $Y$  and $\Phi_{g\circ f}(x)=\Phi_g(f(x))\Phi_f(x)$ $\mu-a.e$ in $Y$.  
\end{proposition}

\demo 

First note that $\mu\circ(g\circ f)$ is a $\sigma$-finite measure in $Y$. 

Now, for each $E\subseteq Y$, $$\int\limits_E\Phi_f(x)\Phi_g(f(x))d\mu=\int\limits_E\Phi_g(f(x))d(\mu\circ f)=\int\limits_{f(E)}\Phi_g(x)d\mu=$$ $$=\int\limits_{f(E)}d(\mu\circ g)=\int\limits_Ed(\mu\circ g\circ f).$$ The first and the third equality are a consequence of the Radon-Nikodym derivative. The other two follow by the change of variable theorem. 
So for each $E\subseteq Y$, $$\int\limits_E\Phi_f(x)\Phi_g(f(x))d\mu=\int\limits_E d(\mu\circ g\circ f)=(\mu\circ g \circ f)(E).$$ So, if $\mu(E)=0$ then $(\mu\circ g\circ f)(E)=0$. By \cite{hewitt}
there exists the Radon-Nikodym derivative $\Phi_{g\circ f}$ of $\mu\circ f\circ g$ with respect to $\mu$ in $Y$ and the equality $(\mu\circ f \circ g)(E)=\int\limits_E\Phi_{g\circ f}(x)d\mu$ holds, for each $E\subseteq Y$. So, for each $E\subseteq Y$, 
$$\int\limits_E\Phi_f(x)\Phi_g(f(x))d\mu=\int\limits_E\Phi_{g\circ f}(x)d\mu,$$ and therefore,$\Phi_f(x)\Phi_g(f(x))=\Phi_{g\circ f}(x)$ $\mu-a.e.$ 
\fim

\section{Representations of Cuntz-Krieger \newline algebras for infinite matrices.}

Representations of the Cuntz-Krieger algebras are of great importance, having applications both to operator algebras and to dynamical systems. In this section we show that for each $A_\infty$-branching system, there exists a representation of the unital Cuntz-Krieger C*-algebra $O_A$ on $\mathcal{B}(L_2(X,\mu))$, the bounded operators on $L_2(X,\mu)$.  

Following \cite{exellaca}, recall that the unital Cuntz-Krieger algebra of an infinite matrix $A$, with $A(i,j)\in \{0,1\}$ and $(i,j)\in \N\times \N$ is the unital universal C*-algebra generated by a family $\{S_i\}_{i\in \N}$ of partial isometries that satisfy:
\begin{enumerate}
\item $S_iS_i^*S_jS_j^*=0$ if $i\neq j$;
\item $S_i^*S_i$ and $S_j^*S_j$ commute, for all $i,j$;
\item $S_i^*S_iS_jS_j^*=A(i,j)S_jS_j^*$, for all $i,j$;
\item $\prod\limits_{u\in U}S_uS_u^*\prod\limits_{v\in V}(1-S_vS_v^*)=\sum\limits_{j=1}^\infty A(U,V,j)S_jS_j^*$, for each pair of finite subsets $U,V\subseteq \N$ such that $A(U,V,j):=\prod\limits_{u\in U}A(u,j)\prod\limits_{v\in V}(1-A(v,j))$ vanishes for all but a finite number of $j's$.

\end{enumerate}   

\begin{theorem}
For a given $A_\infty$-branching system (see \ref{brancsystem}), there exist a *-homomorphism $\pi:O_A\rightarrow \mathcal{B}(L_2(X,\mu))$ such that $\pi(S_i)\phi=\chi_{R_i}\cdot(\Phi_{f_i^{-1}})^{\frac{1}{2}}\cdot\phi\circ F$ for each $\phi\in L_2(X,\mu)$.  
\end{theorem}

\demo

First notice that for a given $\phi\in L_2(X,\mu)$ we have that $$\int\limits_X|\chi_{R_i}(x)\Phi_{f_i^{-1}}(x)^{\frac{1}{2}}\phi(F(x))|^2d\mu=\int\limits_{R_i}\Phi_{f_i^{-1}}(x)|\phi(f_i^{-1}(x))|^2d\mu=$$ $$=\int\limits_{R_i}|\phi(f_i^{-1}(x))|^2d(\mu\circ f_i^{-1})=\int\limits_{D_i}|\phi(x)|^2d\mu\leq \int\limits_{X}|\phi(x)|^2d\mu.$$ To obtain the second equality we have considered the Radon-Nikodym derivative of $\mu\circ f_i^{-1}$ with respect to $\mu$ in $R_ i$ and the last equality is an application of the change of variable theorem.
 
So, we define the operator $\pi(S_i):\Lh(L_2(X,\mu))\rightarrow \Lh(L_2(X,\mu))$ by  $$\pi(S_i)\phi=\chi_{R_i}\cdot (\Phi_{f_i^{-1}})^{\frac{1}{2}}\cdot (\phi\circ F),$$ for each $\phi\in L_2(X,\mu)$. By using the above computation, we see that $\pi(S_i)\in \mathcal{B}(L_2(X,\mu))$. 

Our aim is to show that $\{\phi(S_i)\}_{i\in \N}$ satisfies the relations 1-4 which define the Cuntz-Krieger algebra $O_A$. With this in mind, let us first determine the operator $\psi(S_i)^*$.

For each $\phi,\psi \in L_2(X,\mu)$, 

$$\left\langle\pi(S_i)\phi,\psi \right\rangle=\int\limits_X\chi_{R_i}(x)\Phi_{f_i^{-1}}(x)^{\frac{1}{2}}\phi(F(x))\overline{\psi(x)}d\mu=\int\limits_{R_i}\Phi_{f_i^{-1}}(x)^\frac{1}{2}\phi(f_i^{-1}(x))\overline{\psi(x)}d\mu=...$$  
...by using the change of variable theorem...

$$...=\int\limits_{D_i}\Phi_{f_i^{-1}}(f_i(x))^{\frac{1}{2}}\phi(x)\overline{\psi(f_i(x))}d(\mu\circ f_i)=...$$

...considering the Radon derivative $\Phi_{f_i}$ of $\mu\circ f_i$...

$$...=\int\limits_{D_i}\Phi_{f_i}(x)\Phi_{f_i^{-1}}(f_i(x))^{\frac{1}{2}}\phi(x)\overline{\psi(f_i(x))}d\mu=...$$
...by proposition \ref{prop1}...
$$...=\int\limits_{D_i}\Phi_{f_i}(x)^{\frac{1}{2}}\phi(x)\overline{\psi(f_i(x))}d\mu=\int\limits_{X}\phi(x)\Phi_{f_i}(x)^{\frac{1}{2}}\overline{\psi(f_i(x))}d\mu=\left \langle \phi , \chi_{D_i}\cdot\Phi_{f_i}^{\frac{1}{2}}\cdot(\psi\circ f_i) \right\rangle.$$ Then $$\pi(S_i)^*\psi=\chi_{D_i}\cdot\Phi_{f_i}^{\frac{1}{2}}\cdot(\psi\circ f_i).$$

Using proposition \ref{prop1} again, it is easy to show that 
$$\pi(s_i)^*\pi(S_i)\psi=\chi_{D_i}\cdot \psi =M_{\chi_{D_i}}(\psi)$$ for each $\psi\in L_2(X,\mu)$ (that is, $\pi(S_i)^*\pi(S_i)$ is the multiplication operator by $\chi_{D_i}$). In the same way $\pi(s_i)\pi(S_i)^*=M_{\chi_{R_i}}$.

Now we verify if $\{\pi(S_i)\}_{i\in \N}$ satisfies the relations 1-4, which define the C*-algebra $O_A$.
The first relation follows from the fact that $\mu(R_i\cap R_j)=0$ for $i\neq j$. The second one is trivial.

To see that the third relation is also satisfied, recall that if $A(i,j)=0$ then $\mu(R_j\cap D_i)=0$ and hence $$\pi(S_i)^*\pi(S_i)\pi(S_j)\pi(S_j)^*=M_{\chi_{D_i}}M_{\chi_{R_j}}=M_{\chi_{D_i\cap R_j}}=0,$$ and if $A(i,j)=1$ then $\mu(R_j\setminus D_i)=0$ and hence $$\pi(S_i)^*\pi(S_i)\pi(S_j)\pi(S_j)^*=M_{\chi_{D_i}}M_{\chi_{R_j}}=M_{\chi_{D_i\cap R_j}}=M_{\chi_{R_j}}=\pi(S_j)\pi(S_j)^*.$$ So, for each $i,j\in \N$ 
$$\pi(S_i)^*\pi(S_i)\pi(S_j)\pi(S_j)^*=A(i,j)\pi(S_j)\pi(S_j)^*.$$

To verify the last relation, let $U,V$ be finite subsets of $\N$ such that $A(U,V,j)=1$ only for finitely many $j's$. 

Then, by definition \ref{brancsystem}:5, 
$$M_{\chi_{\left(\bigcap\limits_{u\in U}D_u\bigcap\limits_{v\in V}(X\setminus D_v)\right)}}=M_{\chi_{\left(\bigcup\limits_{A(U,V,j)=1}R_j\right)}}.$$ Note that $$M_{\chi_{\left(\bigcap\limits_{u\in U}D_u\bigcap\limits_{v\in V}(X\setminus D_v)\right)}}=\prod\limits_{u\in U}M_{\chi_{D_u}}\prod\limits_{v\in V}(Id-M_{\chi{D_v}})= $$
$$=\prod\limits_{u\in U}\pi(S_u)^*\pi(S_u)\prod\limits_{v\in V}(Id-\pi(S_v)^*\pi(S_v)).$$ On the other hand,  

$$M_{\chi_{\left(\bigcup\limits_{j\in \N:A(U,V,j)=1}R_j\right)}}=\sum\limits_{j\in \N:A(U,V,j)=1}M_{\chi_{R_j}}=\sum\limits_{j\in \N:A(U,V,j)=1}\pi(S_j)\pi(S_j)^*.$$ This shows that the last relation defining $O_A$ is also verified. 

So, there exist a *-homomorphism $\pi:O_A\rightarrow \mathcal{B}(L_2(X,\mu))$ satisfying   $\pi(S_i)\phi=\chi_{R_i}\cdot(\Phi_{f_i^{-1}})^{\frac{1}{2}}\cdot\phi\circ F$.
\fim

The previous theorem applies only if an $A_\infty$-branching system is given. Our next step is to guarantee the existence of $A_\infty$-branching systems for any matrix $A$. 
First we prove a lemma, which will be helpful in some situations.

\begin{lema}
Let $A$ be a infinite matrix with entries in $\N\times \N$ having no identically zero rows, $(X,\mu)$ be a measure space, and let $\{R_j\}_{j=1}^\infty$ and $\{D_j\}_{j=1}^\infty$ be families of measurable subsets of $X$ such that
\begin{enumerate}
\item [a)] $\mu(R_i\cap R_j)=0$ for all $i\neq j$;
\item [b)] $X\stackrel{\mu-a.e.}{=}\bigcup\limits_{j=1}^\infty R_j$;
\item [c)] $D_i\stackrel{\mu-a.e.}{=}\bigcup\limits_{j\in \N:A(i,j)=1}R_j$;
\end{enumerate}
  Then conditions 4 and 5 of \ref{brancsystem} are satisfied. 
\end{lema}

\demo
Condition 4 follows from $a)$ and $b)$. To show 5 firs we note that $X\setminus D_v\stackrel{\mu-a.e.}{=}\bigcup\limits_{j\in \N:A(v,j)=0}$.
Then, given $U,V$ finite subsets of $N$, we have that
$$\bigcap\limits_{u\in U}D_u\cap \bigcap\limits_{v\in V}(X\setminus D_v)\stackrel{\mu-a.e}{=}\left(\bigcup\limits_{j\in\N:A(u,j)=1 \forall u\in U}R_j\right)\cap \left(\bigcup\limits_{j\in \N:A(v,j)=0 \forall v\in V}R_j \right)\stackrel{\mu-a.e}{=}$$

$$\left(\bigcup\limits_{j\in\N:\prod\limits_{u\in U}A(u,j)=1}R_j\right)\cap \left(\bigcup\limits_{j\in \N:\prod\limits_{v\in V}(1-A(v,j))=1}R_j \right)\stackrel{\mu-a.e}{=}\bigcup\limits_{j\in \N:A(U,V,j)=1}R_j.$$
\fim

\begin{theorem}
For each infinite matrix $A$, without identically zero rows, there \hfill exists \hfill an \hfill  $A_\infty$-branching \hfill system \hfill in \hfill the \hfill measure \hfill space \newline $([0,\infty),\mu)$, where $\mu$ is the Lebesgue measure.
\end{theorem}

\demo
Consider $[0,\infty)$ with the Lebesgue measure $\mu$. Define $R_i=[i,i+1]$ and $D_i=\bigcup\limits_{j:A(i,j)=1}R_j$. Note that $\mu(R_i\cap R_j)=0$ for $i\neq j$. Then, by the previous lemma, conditions 4 and 5 of definition $\ref{brancsystem}$ are satisfied. So, it remains to define maps $f_i:D_i\rightarrow R_i$ and $F:[0,+\infty)\rightarrow [0,+\infty)$ satisfying the conditions of definition $\ref{brancsystem}$. For a fixed $i_0\in \N$ we define $f_{i_0}$ as follows. 
First divide the interval $\stackrel{\circ}{R_{i_0}}$ (where $\stackrel{\circ}{R_{i_0}}$ denotes the interior of $R_{i_0}$) in $\#\{j:A(i_0,j)=1\}$ intervals $I_j$. Then, define $\tilde{f_{i_0}}:\bigcup\limits_{j:A(i_0,j)=1}\stackrel{\circ}{R_j}\rightarrow \bigcup\limits_{j:A(i_0,j)=1}\stackrel{\circ}{I_j}$ such that $\tilde{f_{i_0}}:\stackrel{\circ}{R_j}\rightarrow \stackrel{\circ}{I_j}$ is a $C^1$-diffeomorphism. We now define $f_{i_0}:D_{i_0}\rightarrow R_{i_0}$ by $$f_{i_0}(x)=\left\{\begin{array}{cc}\tilde{f_{i_0}}(x) & \text{ if } x\in \bigcup\limits_{j:A(i_0,j)=1}\stackrel{\circ}{R_j}\\
i_0 & \text{ if } x\in D_{i_0}\setminus \bigcup\limits_{j:A(i_0,j)=1}\stackrel{\circ}{R_j}
\end{array}\right. $$ and
$F:[0,\infty)\rightarrow [0,\infty)$ by 
 $$F(x)=\left\{\begin{array}{cc}\tilde{f_{i_0}}^{-1}(x) & \text{ if } x\in \bigcup\limits_{j:A(i_0,j)=1}\stackrel{\circ}{I_j}\\
0 & \text{ if } x\in R_{i_0}\setminus \bigcup\limits_{j:A(i_0,j)=1}\stackrel{\circ}{I_j}.
\end{array}\right. $$

Note that $f_i$ and $F$ are measurable maps. Moreover, $\mu\circ f_i$ and $\mu\circ f_i^{-1}$ are $\sigma$-finite measures in $D_i$ and $R_i$. Next we show that there exists the Radon-Nikodym derivatives $\Phi_{f_i}$ of $\mu\circ f_i$ with respect do $\mu$ in $D_i$. Let $E\subseteq D_i$ be such that $\mu(E)=0$. To show that $\mu\circ f_i(E)=0$ it is enough to show that $\mu\circ f_i(E\cap(\bigcup\limits_{j:A(i,j)=1}\stackrel{\circ}{R_j}))=0$, and this equality is true by \cite{fernandez}.
 Then, by \cite{hewitt},
there exist the desired nonnegative Radon-Nikodym derivative $\Phi_{f_i}$. 
In the same way there exists the (nonnegative) Radon-Nikodym derivative $\Phi_{f_i^{-1}}$ of $\mu\circ f_i^{-1}$ with respect to $\mu$ in $R_i$. 
We still need to show that $F$ is nonsingular. For this, let $A\subseteq [0,\infty)$ be such that $\mu(A)=0$. Notice that it is enough to prove that $\mu(F^{-1}(A)\cap R_j)=0$ for each $j$. Now $\mu(F^{-1}(A)\cap R_j)=\mu(f_j(A\cap D_j))=0$, (where the last equality follows from the fact that $\mu\circ f_j\ll\mu$ in $D_j$), and hence $\mu(F^{-1}(A))=0$ as desired.

\fim

\begin{corolario} Given an infinite matrix $A$, there exists a representation of $O_A$ in $L_2([0,\infty),\mu)$ where $\mu$ is the Lebesgue measure. If $A$ is $N\times N$ then there exists a representation of $O_A$ in $L_2([0,N),\mu)$ where $\mu$ is the Lebesgue measure.
\end{corolario}

\section{The Perron-Frobenious Operator}

We now describe the Perron-Frobenious operator using the representations introduced in the previous section.

\begin{theorem}
Let $(X,\mu)$ be a measure space with a branching system as in definition $\ref{brancsystem}$ and let $\varphi\in L_1(X,\mu)$ be such that $\varphi(x)\geq 0$ $\mu$-a.e.. 
\begin{enumerate}
\item If $supp(\varphi)\subseteq \bigcup\limits_{i=1}^NR_j$, then $$P_F(\varphi)=\sum\limits_{i=1}^N\left(\pi(S_i^*)\sqrt{\varphi}\right)^2.$$

\item If $supp(\varphi)\subseteq \bigcup\limits_{i=1}^\infty R_j$, then $$P_F(\varphi)=\lim_{N\rightarrow \infty}\sum\limits_{i=1}^N\left(\pi(S_i^*)\sqrt{\varphi}\right)^2,$$ where the convergence occurs in the norm of $L_1(X,\mu)$.
\end{enumerate}
\end{theorem}

\demo The first assertion follow from the fact that for each measurable set $A\subseteq X$, $\int\limits_AP_F(\varphi)(x)d\mu=\int\limits_A\sum\limits_{i=1}^N\left(\pi(S_i^*)\sqrt{\varphi}(x)\right)^2d\mu$. To prove this equality, we will use the Radon-Nikodym derivative of $\mu\circ f_i$, the change of variable theorem and the fact that $F^{-1}(A)\cap R_i=f_i(A\cap D_i)$. Given $A\subseteq X$ a measurable set we have that

$$\sum\limits_{i=1}^N\int\limits_A\left(\pi(S_i^*)\sqrt{\varphi}(x)\right)^2d\mu=\sum\limits_{i=1}^N\int\limits_A\chi_{D_i}(x)\Phi_{f_i}(x)\varphi(f_i(x))d\mu=$$
$$=\sum\limits_{i=1}^N\int\limits_{A\cap D_i}\Phi_{f_i}(x)\varphi(f_i(x))d\mu=\sum\limits_{i=1}^N\int\limits_{A\cap D_i}\varphi(f_i(x))d(\mu\circ f_i)=\sum\limits_{i=1}^N\int\limits_{f_i(A\cap D_i)}\varphi(x)d\mu=$$ 
$$=\sum\limits_{i=1}^N\int\limits_{F^{-1}(A)\cap R_i}\varphi(x)d\mu=\sum\limits_{i=1}^N\int\limits_{F^{-1}(A)}\chi_{R_i}\varphi(x)d\mu=\int\limits_{F^{-1}(A)}\sum\limits_{i=1}^N\chi_{R_i}\varphi(x)d\mu=$$
$$=\int\limits_{F^{-1}(A)}\varphi(x)d\mu=\int\limits_AP_F(\varphi)(x)d\mu.$$

We now prove the second assertion. For each $N\in \N$, define $\varphi_N:=\sum\limits_{i=1}^N\chi_{R_i}\cdot \varphi.$ 
Note that $(\varphi_N)_{N\in \N}$ is an increasing sequence, bounded above  by $\varphi$. 
Then 
$$\lim\limits_{N\rightarrow\infty}\int\limits_XP_F(\varphi_N)(x)d\mu=\lim\limits_{N\rightarrow \infty}\int\limits_X\varphi_N(x)d\mu=...$$
...by the Lebesgue's Dominated Convergence Theorem...

$$=\int\limits_X\varphi(X)d\mu=\int\limits_XP_F(\varphi)(x)d\mu.$$

Moreover, the sequence $(P_F(\varphi_N))_{N\in \N}$ is $\mu$ - a. e. increasing and bounded above by $P_F(\varphi)$.  

Then, $$\lim\limits_{N\rightarrow\infty}\|P_F(\varphi)-P_F(\varphi_N)\|_1=\lim\limits_{N\rightarrow\infty}\int\limits_X|P_F(\varphi)(x)-P_F(\varphi_N)(x)|d\mu=$$ 
$$=\lim\limits_{N\rightarrow\infty}\int\limits_XP_F(\varphi)(x)-P_F(\varphi_N)(x)d\mu=0.$$
Therefore, $\lim\limits_{N\rightarrow \infty}P_F(\varphi_N)=P_F(\varphi)$. By the first assertion, 
$P_F(\varphi_N)=\sum\limits_{i=1}^N\left(\pi(S_i^*)\sqrt{\varphi_N}\right)^2$, and a simple calculation shows that
$$\sum\limits_{i=1}^N\left(\pi(S_i^*)\sqrt{\varphi_N}\right)^2=\sum\limits_{i=1}^N\left(\pi(S_i^*)\sqrt{\varphi}\right)^2.$$

So, we conclude that $$\lim\limits_{N\rightarrow \infty}\sum\limits_{i=1}^N\left(\pi(S_i^*)\sqrt{\varphi}\right)^2=P_F(\varphi).$$
\fim

\begin{theorem} Let $A$ be a matrix such that each row has a finite number of 1's and let $(X,\mu)$ be an $A_\infty$-branching system. Suppose $\mu(R_i)< \infty$ for each $i$ (so that $\chi_{R_i}\in L_1(X,\mu)$). Moreover, suppose $\Phi_{f_i}$ is a constant positive function for each $i$, say $\Phi_{f_i}=b_i$ (for example, if $f_i$ is linear). Let $W\subseteq L_1(X,\mu)$ be the vector subspace $$W=span\{\chi_{R_i}: i\in \N\},$$ that is, $W$ is the subspace of all finite linear combinations of $\chi_{R_i}$. Then the Perron-Frobenius operator restricted to $W$, ${P_F}_{|_W}:W\rightarrow W$, has a matrix representation given by $A^T B$, where $B$ is the diagonal infinite matrix with nonzero entries $B_{i,i}=b_i$. 
\end{theorem}

Although $A$ and $B$ are infinite matrices, we are considering the matrix multiplication $A^TB$ as the usual multiplication for finite matrices, since $B$ is column-finite. 

\demo Since each row $z$ of $A$ has a finite number of 1's, then, by definition \ref{brancsystem}:5, taking $Z=\{z\}$ and $Y=\emptyset$, we obtain $D_z\stackrel{\mu-a.e}{=}\bigcup\limits_{j:A(z,j)=1}R_j$
so that $\chi_{D_z}=\sum\limits_{j:A(z,j)=1}\chi_{R_j}$. Note that 

$$P_F(\chi_{R_z})=b_z\chi_{D_z}=\sum\limits_{j:A(z,j)=1}b_zR_j,$$ and so the element $(j,z)$ of the matrix representation of ${P_F}_{|_W}$ is $b_zA(z,j)$. \fim

\section{Examples}
\begin{example} $O_\infty$ ($O_A$ where all entries of the matrix $A$ are 1).\end{example} Consider $X=[0,1]$ with Lebesgue measure and define $D_i=[0,1]$, for $i=1,2,\ldots$. To define the $R_i$´s we first need to define recursively the following sequences in $X$: Let $a_1=0$, $a_i=a_{i-1}+\frac{1}{2^i}$, $i=2,3,\ldots$ and let $b_i=\frac{a_i+a_{i+1}}{2}$, $i=1,2,\ldots$. Now define $R_i=[a_{\frac{i+1}{2}},b_{\frac{i+1}{2}}]$ for $i$ odd and $R_i=[b_{\frac{i}{2}}, a_{\frac{i}{2}+1}]$ for $i$ even and define a map $F$ on $X$ by $F(x)=\frac{x}{b_{\frac{i+1}{2}}-a_{\frac{i+1}{2}}}+\frac{a_{\frac{i+1}{2}}}{a_{\frac{i+1}{2}}-b_{\frac{i+1}{2}}}$ for $x\in R_i$, $i$ odd and $F(x)=\frac{x}{a_{\frac{i}{2}+1}-b_{\frac{i}{2}}}+\frac{b_{\frac{i}{2}}}{b_{\frac{i}{2}}-a_{\frac{i}{2}+1}}$ for $x\in R_i$, $i$ even. Notice that $F$ is nothing more than an affine transformation that takes the interval $R_i$ onto $D_i=[0,1]$, as shown in the picture below:

\begin{center}
\epsfxsize=7cm    
\epsfbox{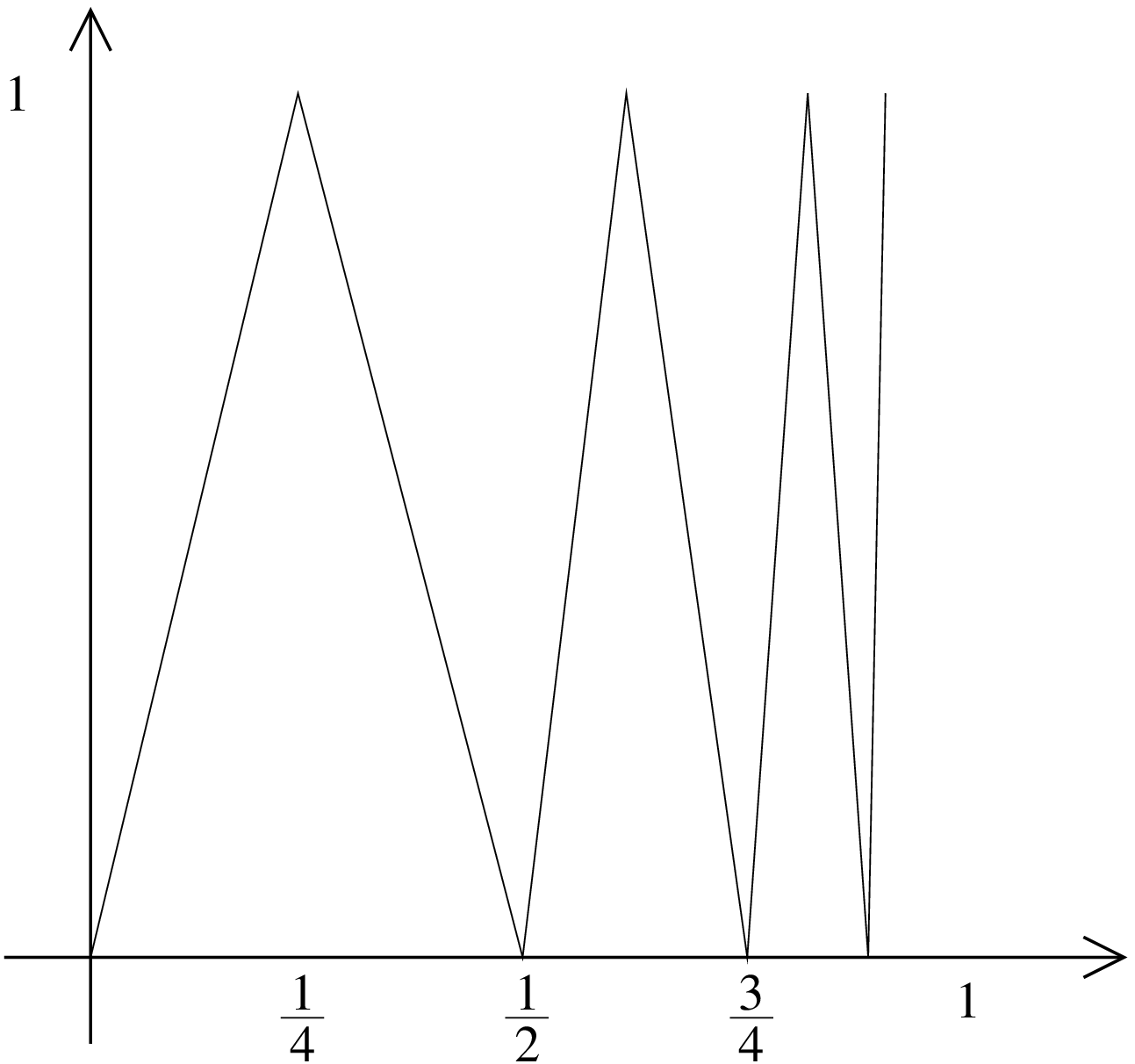}
\end{center}

Finally, let $f_i=(F_{|_{R_i}})^{-1}$. Then $(\{f_i\}_{i=1}^\infty, \{D_i\}_{i=1}^\infty)$ is an $A_\infty$ branching system and hence induces a representation of the Cuntz-Krieger algebra $O_\infty$.

\begin{example} 
\end{example}

Let $X$ be the measure space $[0,\infty)$, with the Lebesgue measure. 
Consider the map $F:[0,\infty)\rightarrow [0,\infty)$ defined by  $F(x)=\frac{i}{2}(x-i)^2$ for $x\in [i-1,i]$ and $i$ odd and $F(x)=[\frac{i}{2}](x-(i-1))^2$ for $x\in [i-1,i]$ and $i$ even ( $[\frac{i}{2}]$ is the least integer greater than or equal to $\frac{i}{2}$). Below we see the graph of $F$.

\begin{center}
\epsfxsize=11cm    
\epsfbox{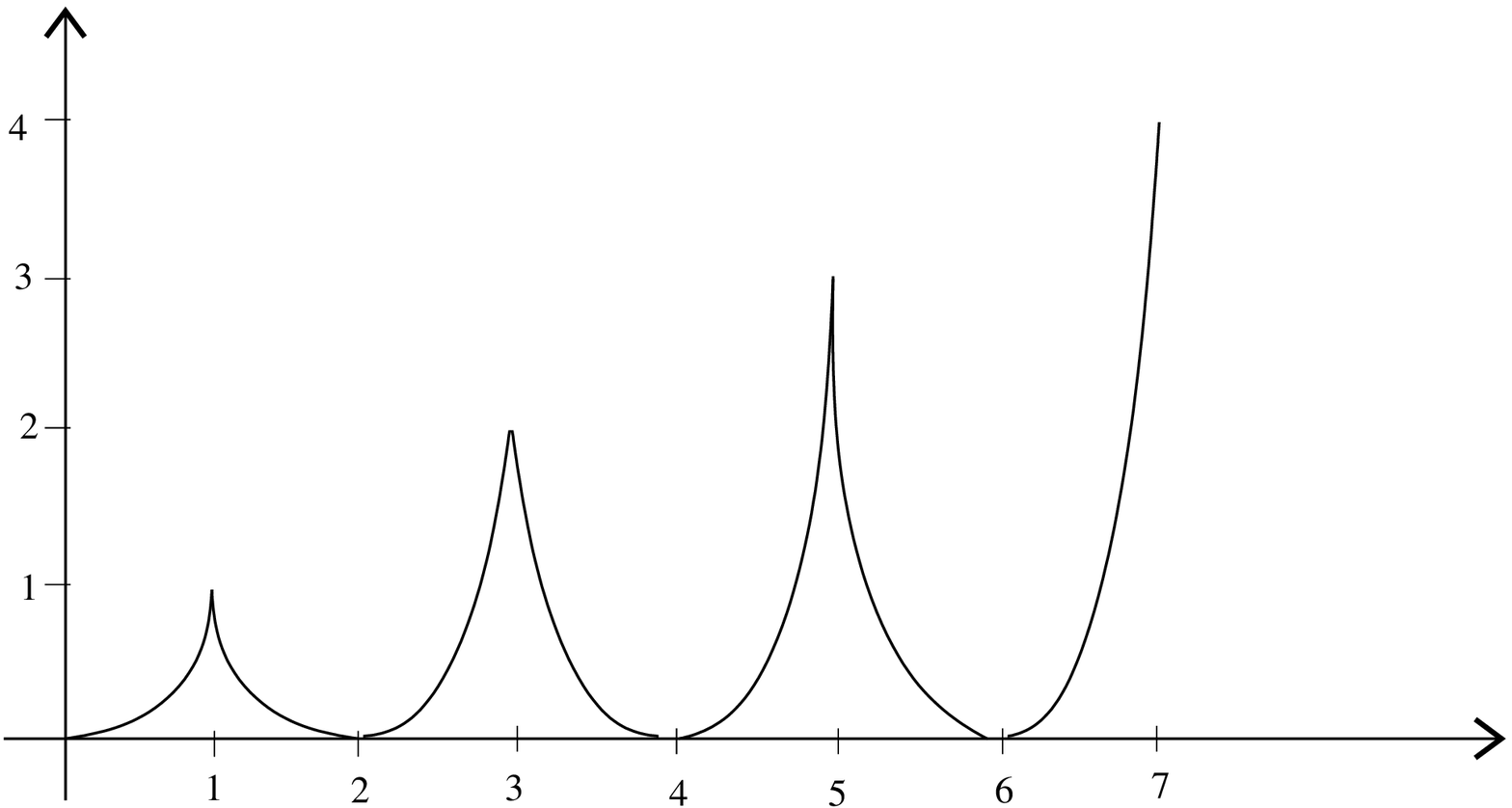}
\end{center}

Define $R_i=[i-1,i]$ for $i=1,2,3,...$, set $D_i=[0,[\frac{i}{2}]]$ and let $f_i:D_i\rightarrow R_i$ be defined by $f_i=(F_{|_{R_i}})^{-1}$. Then $(\{f_i\}_{i=1}^\infty,\{D_i\}_{i=1}^\infty)$ is an $A_\infty$ branching system. This branching system induces a representation of the C*-algebra $O_A$, for 
$$A=\left(\begin{array}{ccccc}
1 & 0 & 0 & 0 & \cdots\\
1 & 0 & 0 & 0 & \cdots\\
1 & 1 & 0 & 0 & \cdots \\
1 & 1 & 0 & 0 & \cdots \\
1 & 1 & 1 & 0 & \cdots \\
\vdots & \vdots & \vdots & 
\end{array}\right).$$

\vspace{1.5pc}

D. Goncalves, Departamento de Matemática, Universidade Federal de Santa Catarina, Florianópolis, 88040-900, Brasil

Email: daemig@gmail.com

\vspace{0.5pc}
D. Royer, Departamento de Matemática, Universidade Federal de Santa Catarina, Florianópolis, 88040-900, Brasil

Email: royer@mtm.ufsc.br
\vspace{0.5pc}

\end{document}